\documentclass[12pt,leqno]{article}
\usepackage{amsmath, amsthm, amsfonts, amssymb, color}
\usepackage{mathrsfs}
\theoremstyle{definition}

\newcommand{\scr}[1]{\mathscr #1}
\definecolor{wco}{rgb}{0.5,0.2,0.3}

\numberwithin{equation}{section} \theoremstyle{remark}

\newcommand{\ua}{\uparrow}

\title{{\bf Stochastic Generalized Porous Media Equations with Reflection}\footnote{Supported in
 part by NNSFC(11131003), SRFDP, the Fundamental Research Funds for the Central Universities and the DFG through IRTG 1132.} }
\author{
{\bf  Michael R\"ockner$^{b)}$, Feng-Yu Wang$^{a),d)}$\footnote{Corresponding author.
wangfy@bnu.edu.cn; F.Y.Wang@swansea.ac.uk}, Tusheng Zhang$^{c)}$}\\
\footnotesize{$^{a)}$School of Mathematical Sciences,
Beijing Normal
University, Beijing 100875, China}\\
\footnotesize{$^{b)}$Department of Mathematics, Bielefeld
University, D-33501 Bielefeld, Germany}\\
\footnotesize{$^{c)}$\footnotesize{c) School of Mathematics, University of Manchester, Oxford Road, Manchester M13 9PL, UK}}\\
\footnotesize{$^{d)}$Department of Mathematics, Swansea University,
Singleton Park, SA2 8PP, UK} }
\begin{document}
\def\R{\mathbb R}  \def\ff{\frac} \def\ss{\sqrt} \def\B{\mathbb
B} \def\BB{\scr B} \def\LL{\Lambda}
\def\N{\mathbb N} \def\kk{\kappa} \def\m{{\bf m}}
\def\dd{\delta} \def\DD{\Delta} \def\vv{\varepsilon} \def\rr{\rho}
\def\<{\langle} \def\>{\rangle} \def\GG{\Gamma} \def\gg{\gamma}
  \def\nn{\nabla} \def\pp{\partial} \def\tt{\tilde}
\def\d{\text{\rm{d}}} \def\bb{\beta} \def\aa{\alpha} \def\D{\scr D}
\def\EE{\scr E} \def\si{\sigma} \def\ess{\text{\rm{ess}}}
\def\beg{\begin} \def\beq{\begin{equation}}  \def\F{\scr F}
\def\Ric{\text{\rm{Ric}}} \def\Hess{\text{\rm{Hess}}}
\def\e{\text{\rm{e}}} \def\ua{\underline a} \def\OO{\Omega}  \def\oo{\omega}
 \def\tt{\tilde} \def\Ric{\text{\rm{Ric}}}
\def\cut{\text{\rm{cut}}} \def\P{\mathbb P} \def\ifn{I_n(f^{\bigotimes n})}
\def\C{\scr C}      \def\aaa{\mathbf{r}}     \def\r{r}
\def\gap{\text{\rm{gap}}} \def\prr{\pi_{{\bf m},\varrho}}  \def\r{\mathbf r}
\def\Z{\mathbb Z} \def\vrr{\varrho} \def\ll{\lambda}
\def\L{\scr L}\def\Tt{\tt} \def\TT{\tt}\def\II{\mathbb I}
\def\i{{\rm i}}\def\Sect{{\rm Sect}}\def\E{\mathbb E} \def\H{\mathbb H}
\def\M{\scr M}\def\Q{\mathbb Q} \def\texto{\text{o}} \def\i{{\rm i}}
\def\O{\scr O}
\maketitle
\begin{abstract} A non-negative Markovian solution is constructed for a class of stochastic generalized porous media equations with reflection. To this end,
some regularity properties and a comparison theorem are proved for stochastic generalized porous media equations, which are   interesting by themselves.   Invariant probability measures and ergodicity of the solution are also investigated.
\end{abstract} \noindent
 AMS subject Classification:\ 60J75, 47D07.   \\
\noindent
 Keywords:     Stochastic porous media equation, reflection, regular solution, comparison theorem.

 \vskip 2cm

\section{Introduction}

Let $E$ be a locally compact separable metric space with Borel $\si$-field $\BB$ and let $\mu$ be a probability measure on $(E,\BB)$. Let $(L,\D(L))$ be a symmetric Dirichlet operator on $L^2(\mu)$ with empty essential spectrum and regular  Dirichlet form  $(\EE,\D(\EE)).$  Note that  we may allow the Dirichlet form to be merely quasi-regular by ``local compactification" from the book \cite{MR}. Let
$\{\ll_i\}_{i\ge 1}$ be all eigenvalues of $-L$ counting multiplicities in increasing order such that $\ll_1>0$, and let $\{e_i\}_{i\ge 1}$ be the corresponding normalized eigenfunctions.

In this paper we investigate the stochastic generalized porous medium equation with reflection of the type
\beq\label{1.1} \d X_t =\{L\Psi(X_t)+\Phi(X_t)\}\d t +\si_t\d W_t +\d \eta_t,\ \ t\in [0,T],\end{equation} where $T>0$ is a fixed constant, $W_t$ is a cylindrical Brownian motion on $L^2(\mu)$ with respect to a complete filtered probability space $(\OO,\F,\{\F_t\}_{t\ge 0},\P)$, $\Psi,\Phi\in C(\R)$, and $\sigma_t\in L^2(\OO\times [0,T]\to \scr L_{HS})$
 is progressively measurable. Here
$\scr L_{HS}$ is the space of all Hilbert-Schmidt linear operators on $L^2(\mu)$. We will use $\<\cdot,\cdot\>$ and $\|\cdot\|$ to denote the inner product and the norm in $L^2(\mu)$, and denote  the norm in $L^p(\mu)$ for $p\ge 1$ by $\|\cdot\|_p$.    Moreover, let $H^1=\D(\EE)$ with inner product
$\<u,v\>_{H^{1}}= \EE(u,v)$, and let $H^{-1}$ be the dual space of $H^1$ w.r.t. $L^2(\mu)$. For simplicity we also use $\<\cdot,\cdot\>$ for the dualization $_{H^{-1}}\<\cdot,\cdot\>_{H^1}$ between $H^{-1}$ and $H^1$. Since the essential spectrum of $L$ is empty, we have $\ll_i\uparrow \infty$ as $i\uparrow\infty$ and thus,  the embedding $H^1\subset L^2(\mu)$ is compact.   As usual $L$ extends  to an operator $L: H^1\to H^{-1},$ again denoted by $L$, and below by $L$ we always mean this extension.

\vskip 0.3cm
For stochastic partial differential equations with reflection driven by space-time white noise, we refer the reader to
 \cite{DP1}, \cite{NP}, \cite{XZ}, \cite{ZA} and \cite{Z}. The situation here is drastically different from that in the above references because of the  presence of the the non-linear operator $L\Psi$, in (\ref{1.1}).
\vskip 0.3cm

The motivation to study reflection problems of type (\ref{1.1}) is that eventually we would like to extend our results from this paper to stochastic fast diffusion equations, where $\Psi(s)= |s|^{r-1}s$ with $r\in (0,1)$, or to so-called ``self-organized criticality" models, where $\Psi$ is a Heaviside function or a product of this with the identity. Such type of singular stochastic porous media equations have intensively been studied in \cite{0,01,02} and \cite{7'}, proving in addition to well-posedness that extinction occurs with strictly positive probability. However, in contrast to the latter papers, instead of linear multiplicative noise, we would  like to study the class of additive noise. This was suggested in \cite{03} by A. Diaz-Guilera, who derived equations of type (\ref{1.1}) with $\eta\equiv 0$ as models for the phenomenon of self-organized criticality, where $X_t, t\ge 0$, has the interpretation of energy. However, for all types of additive noise suggested in \cite{03} the solution $X_t, t\ge 0,$ can take (depending on $\oo\in \OO$) arbitrarily negative values, which is somehow in contradiction of its interpretation as energy. Our ``penalization" term, however, guaranties nonnegative solutions. Therefore, we suggest our equation   (\ref{1.1}) as a more realistic version of the one suggested in \cite{03}, where it was also pointed out that the case $\Psi(s)= |s|^{r-1}s$ for $r=3$ is an interesting special case, since it is the simplest fulfilling all symmetry restriction suggested by physical considerations. And this case is covered by the main result in this paper. As also pointed out in \cite{03}, this polynomial case is, however, too much of a simplification and quite far from ``self-organized criticality" models as $\Psi(s)=\mathbf H(s)$ or $\Psi(s)= s\mathbf H(s)$ with $\mathbf H$ being the Heaviside function. In contrast to the polynomial case the latter one namely exhibits extinction of solutions as shown in \cite{0,02} and \cite{7'}, while the polynomial case does not. Therefore, in a future paper we plan to extend our results to singular cases as $\Psi(s)= |s|^{r-1} s$ for $r\in (0,1), \Psi(s)=\mathbf H(s),$ or $\Psi(s)= s\mathbf H(s).$

 Let $\scr M_c$ be the space of all locally finite measures on $E$,
equipped with the vague topology induced by $f\in C_0(E)$, where     $C_0(E)$ be the set of all continuous functions on $ E$ with compact support.
\beg{defn} An element $u\in H^{-1}$ is called non-negative,   denoted by $u\ge 0$, if $_{H^1}\<f,u\>_{H^{-1}}\ge 0$ holds for any non-negative $f\in H^1$. For $u_1,u_2\in H^{-1}$, we write $u_1\ge u_2$ if $u_1-u_2\ge 0.$ A process $u_t$ in $H^{-1}$ is called increasing if $u_t\ge u_s$ for $t\ge s.$ \end{defn}

It is easy to verify that $u_1\ge u_2$ and $u_2\ge u_1$ if and only if $u_1=u_2.$ Thus, $H^{-1}$ is a partially ordered space.

\beg{defn} \label{D1.1} A pair $(X,\eta):=(X_t,\eta_t)_{t\ge 0}$ is called a solution to (\ref{1.1}), if \beg{enumerate}
\item[(1)] $X$ is a non-negative,  adapted   process on $L^2(\mu)$, which is c\'adl\'ag   in   $H^{-1}$,  such that for any $T>0,$
$\Psi(X_{\cdot})|_{[0,T]}\in L^2(\OO\times [0,T]\to H^1; \P\times \d t)$  and $\Phi(X_{\cdot})|_{[0,T]}\in L^2(\OO\times [0,T]\to H^{-1}; \P\times\d t)$;
\item[(2)] $\eta=(\eta_t)_{t\in [0,T]}$ is a right-continuous, non-negative, increasing adapted process in $H^{-1}$, which determines a  unique
 adapted process $\nu:=(\nu_t)_{t\ge 0}$  in $\scr M_c$,    right-continuous in the   topology of set-wise convergence,   such that
\beq\label{D1} \int_Ef(z)\nu_t(\d z)= \ _{H^{-1}}\<\eta_t, f\>_{H^1},\ \ f\in C_0(E)\cap H^1,\ \ t\ge 0.\end{equation}
Moreover,  $_{H^{1}}\<\Psi(X_t),\eta_t\>_{H^{-1}}=0,\ \P\times \d t$-a.e.;
\item[(3)]
$\P$-a.s.
$$X_t = X_0 +\int_0^t \Big\{L\Psi(X_s)+\Phi(X_s)\Big\}\d s +\int_0^t \si_s\d W_s+\eta_t,\ \ t\ge 0$$ holds on $H^{-1}.$  \end{enumerate}
\end{defn}


To construct a solution to (\ref{1.1}) using a natural approximation argument, we shall need the following assumptions.
\beg{enumerate} \item[{\bf (A1)}]  For any $f\in C_0(E)$, there exists $\tt f\in H^1\cap C_0(E)$ such that for any $\vv>0$,
$$|f-f_\vv|\le \vv \tt f $$ holds for some $f_\vv\in H^1\cap C_0(E).$
\item[{\bf (A2)}] $\Phi: \R\to\R$ is Lipschitz continuous such that
$$\<\Phi(u)-\Phi(v), (-L)^{-1}(u-v)\>\le c\|u-v\|_{H^{-1}}^2,\ \ u,v\in L^{r+1}(\mu),$$ holds for some constant $c>0$; and  $\Psi\in C^1(\R)$ with $\Psi(0)=0$  and there exist constants $r\ge 1, c_1'\ge 0$ and $c_1,c_2>0$ such that for any $s_1,s_2,s\in\R$,
$$(s_2-s_1)(\Psi(s_2)-\Psi(s_1))\ge  \{c_1|s_2-s_1|^{r+1}+c_1' |s_2-s_1|^2\},\  0\le \Psi'(s) \le c_2(1+|s|^{r-1}).$$
\item[{\bf (A3)}] $\{e_i\}_{i\ge 1}\subset L^{r+1}(\mu)$ and for any $T>0$ there exists a constant $C>0$ such that
$$ \sup_{t\in [0,T]} \|\si_t\|_{\scr L_{HS}}^2+ \sup_{n\ge 1}  \int_{[0,T]\times E} \Big\{\sum_{i=1}^n\Big(\sum_{k=1}^n\<\si_t e_i, e_k\>e_k\Big)^2\Big\}^{(1+r)/2}\d t\d\mu\le C.$$\end{enumerate}

Here is the main result of the paper, where   for a measure $\nu$ on $(E,\scr B)$ and a $\nu$-integrable or nonnegative $\scr B$-measurable function $f: E\to \R$ we set $\nu(f)= \int_Ef\d\nu.$

\vskip 0.3cm
\beg{thm}\label{T1.1} Assume that {\bf (A1)}, {\bf (A2)} and {\bf (A3)} hold. If either $c_1'>0$, or $\Psi(s)= c|s|^{r-1}s$ and $\Phi(s)= c's$ for some constants $c>0$ and $c'\in\R.$ Then: \beg{enumerate}
\item[$(1)$] For any $X_0=x\in L^{1+r}(\OO\to L_+^{1+r}(\mu),\F_0;\P)$, $(\ref{1.1})$ has a  solution $(X(x),\eta(x))$ in the sense of Definition $\ref{D1.1}$ such that
\beq\label{PP1.1}\E \sup_{t\in [0,T]}\|X_t(x)\|^2 + \sup_{t\in [0,T]}\E\|X_t(x)\|_{1+r}^{1+r}<\infty.\end{equation}
If $r=1$ and $\Psi$ is linear, then the solution to $(\ref{1.1})$ is unique.
\item[$(2)$] If $\si_t =\si$ is deterministic and independent of $t$, the family $(X(x))_{x\in L^{1+r}_+(\mu)}$ is time-homogeneous Markovian, i.e. for any $f\in C_b(L^{1+r}(\mu))$ and $0<s<t\le T$,
\beq\label{M}\E\big(f(X_t(x))\big|\si(X_u(x): u\le s)\big)= (P_{t-s}f)(X_s(x)),\end{equation} where
$$P_uf(z):= \E f(X_u(z)),\ \ \ u\in [0,T], z\in L^{1+r}(\mu).$$
\item[$(3)$] Let $\si_t=\si$ be deterministic and independent of $t$, and let  $K\in \R$ be such that
\beq\label{PH} \<\Phi(x)-\Phi(y), 1_{\{x-y>0\}}\>\le K\mu\big((x-y)^+\big),\ \ x,y\in L^{1+r}(\mu),\end{equation}
then  $X_t(x)$ is $L^1$-Lipschitz continuous in $x$, i.e. $\P$-a.s.
\beq\label{PH1} \|X_t(x)-X_t(y)\|_1\le \e^{Kt} \|x-y\|_1,\ \ \ x,y\in L^{1+r}(\mu).\end{equation}\end{enumerate} Consequently, $P_t$ extends to a unique Markov Lipschitz-Feller semigroup on $L^1(\mu)$. If
\beq\label{LLL} \EE(\Psi(x),x)-\<\Phi(x),x\>\ge c_1\EE(x,x)-c_2\end{equation} holds for some constants $c_1,c_2>0$ and all $x\in \D(\EE)$ such that $\Psi(x)\in\D(\EE)$, then
$P_t$   has an invariant probability measure $\pi$ with $\pi(\|\cdot\|_{H^1}^2)<\infty$  and
\beq\label{PH2}|P_t f(x)-\pi(f)|\le {\rm Lip}_1(f)\e^{Kt} \int_{L^1(\mu)}\|x-y\|_1\pi(\d y),\ \ \ x\in L^1(\mu),t\ge 0,\end{equation}  where $f: L^1(\mu)\to \R$ is Lipschitz and ${\rm Lip}_1(f)$ is its Lipschitz constant.   In particular, $P_t$ converges  exponentially fast to $\mu$ if $K<0.$ \end{thm}

To illustrate this result, let us consider the following simple example.

\paragraph{Example 1.1.} Let $E\subset \R^d$ be a bounded open domain, and let $L=\DD$ be the Dirichlet Laplacian on $E$, which has discrete spectrum with $\ll_1>0$. Let $\Phi$ be a Lipschitz continuous function on $\R$ and let $$\Psi(s)= \aa_1s|s|^{r-1}+\aa_2 s +\aa_3 |s|^{r'-1}s$$ for some constants $r\ge 1, r'\in (0,1), \aa_1 >0, $ and $\aa_2,  \aa_3\ge 0.$ Finally, let $\si_t=\si$ be deterministic and independent of $t$ such that
$$\si e_i= q_ie_i,\ \ i\ge 1,$$ where $\{q_i\}_{i\ge 1}\subset [0,\infty)$ such that
\beq\label{F0}\sum_{i=1}^\infty i^{(r^2-1)/(2r+4)}   q_i^2 <\infty.\end{equation} Then $\Psi(0)=0$ and {\bf (A1)}, {\bf (A2)} and {\bf (A3)} hold.

\beg{proof} {\bf (A2)} is trivial for the specific choice of $\Psi$.  Noting that \beq\label{CCC}\|e_i\|_{r+1}\le c \ll_i^{d(r-1)/4(r+1)},\ \ll_i\le ci^{2/d},\ \ \ i\ge 1,\end{equation} holds for some constant $c>0$, we have
\beg{equation*}\beg{split} &\sup_{n\ge 1} \E\int_{[0,T]\times E} \Big\{\sum_{i=1}^n\Big(\sum_{k=1}^n\<\si_t e_i, e_k\>e_k\Big)^2\Big\}^{(1+r)/2}\d t\d\mu   =  T\int_{  E}  \Big\{\sum_{i=1}^\infty q_i^2 e_i^2  \Big\}^{(1+r)/2}\d\mu\\
&\le T \bigg(\sum_{i=1}^\infty
q_i^2\|e_i\|^{1+r}_{1+r}i^{-(r-1)^2/(2r+4)}\bigg)\bigg(\sum_{i=1}^\infty q_i^2i^{(r^2-1)/(2r+4)}\bigg)^{(r-1)/2},\end{split}\end{equation*} which is finite
according to (\ref{F0}) and  (\ref{CCC}).  Therefore, {\bf (A3)} holds.
Finally,  {\bf (A1)} is trivial in the present situation.
\end{proof}

To construct the desired solution to (\ref{1.1}), we first present some preparations in Section 2 concerning regularity properties and a comparison theorem for stochastic generalized porous media equations, which are  interesting   by themselves. A complete proof of Theorem \ref{T1.1} is given in Section 3.

\section{Preparations}

In this section we first consider regularity of solutions to stochastic generalized porous media equations, then prove a comparison theorem and the $L^1$-Lipshitz continuity for the solution. Finally, to construct the ``local time" $\eta$ as a locally finite measure on $[0,T]\times E$, we prove a  suitable new version of the Riesz-Markov representation theorem.  We note that the regularity of solutions for stochastic generalized porous media equations have been investigated  e.g. in \cite{RW08, Gess}
for either linear $\Phi$ or $\Phi=0$. But to approximate the equation with reflection, a non-linear Lipschitzian term $s\mapsto n s^-$ will be  included in $\Phi^{(n)}$ (see Section 3).

\subsection{Regular solution of stochastic generalized porous media equations}

In this subsection we consider the following equation with multiplicative noise:
\beq\label{2.1} \d X_t = \{L\Psi(X_t)+\Phi(X_t)\}\d t +\si_t(X_t)\d W_t,\end{equation} where
$\si: \OO\times [0,\infty)\times L^2(\mu)\to \L_{LS}$ is progressively measurable such that

\beg{enumerate}\item[{\bf (A3')}] $\{e_i\}_{i\ge 1}\subset L^{r+1}(\mu)$ and for any $T>0$ there exists a constant $C>0$ such that for any $u,v\in L^2(\mu)$,
\beg{equation*}\beg{split} &\sup_{t\in [0,T]} \|\si_t(u)-\si_t(v)\|^2_{\scr L_{HS}(L^2(\mu);H^{-1})}\le C\|u-v\|^2_{H^{-1}},\\
 &\sup_{t\in [0,T]} \|\si_t(u)\|_{\scr L_{HS}}^2+ \sup_{n\ge 1}  \int_{[0,T]\times E} \Big\{\sum_{i=1}^n\Big(\sum_{k=1}^n\<\si_t(u) e_i, e_k\>e_k\Big)^2\Big\}^{\ff{1+r}2}\d t\d\mu\le C.\end{split}\end{equation*}\end{enumerate}

\beg{defn}\label{D2.1} A continuous adapted process $X=(X_t)_{t\ge 0}$ on $H^{-1}$ is called  a solution to   (\ref{2.1}), if
$\P$-a.s.
$$X_t= X_0+ \int_0^t \big\{L\Psi(X_s)+\Phi(X_s)\big\}\d s + \int_0^t \si_s(X_s)\d W_s,\ \ \ t\ge 0$$ holds on $H^{-1}$. The solution is called regular, if it is a right-continuous process on $L^2(\m)$ and    for any $T>0$,
$$\E \sup_{t\in [0,T]}\|X_t\|^2 +\sup_{t\in [0,T]}\E\|X_t\|_{1+r}^{1+r} + \E \int_0^T \Big\{\EE(\Psi(X_t),\Psi(X_t))\Big\}\d t<\infty.$$
\end{defn}

\beg{thm}\label{T2.1} Assume {\bf (A1)}, {\bf (A2)} and {\bf (A3')}. For any $\F_0$-measurable random variable  $X_0$ on $H^{-1}$ with $\E \|X_0\|_{H^{-1}}^2<\infty$, the equation $(\ref{2.1})$ has a unique solution such that
\beq\label{ADD} \E\sup_{t\in [0,T]}\|X_t\|_{H^{-1}}^2+\E\int_0^T\|X_t\|_{1+r}^{1+r}\d t <\infty,\ \ \ T>0.\end{equation}
If moreover   $\E\|X_0\|_{1+r}^{1+r}<\infty$, then:
\beg{enumerate} \item[$(1)$] When $c_1'>0$, the solution is regular, continuous in $L^2(\m)$, and  the It\^o formula
\begin{eqnarray}
\|X_t\|^2 &=&\|X_0\|^2 -2\int_0^t \EE(X_s, \Psi(X_s))\d s +2\int_0^t \<X_s, \Phi(X_s)\>\d s \nonumber\\
&& +2\int_0^t \<X_s, \si_s(X_s)\d W_s\>+\int_0^t||\si_s(X_s)||_{\L_{LS}}^2ds,\ \ t\in [0,T]\nonumber
\end{eqnarray}
 holds.
 \item[$(2)$] When $\Psi(s)= c|s|^{r-1}s, \Phi(s)=c's$ for some constants $c>0$ and $c'\in\R$, the solution is regular. \end{enumerate}
\end{thm}

To prove this theorem, we first prove that there exists a  unique solution in $H^{-1}$ using a general result in \cite{RRW07}. Then we show that this solution is indeed regular.

\beg{lem}\label{L2.1} Assume {\bf (A1)}, {\bf (A2)} and {\bf (A3')}. For any $X_0\in L^2(\OO\to H^{-1},\F_0; \P)$, there exists a unique continuous adapted process $X=(X_t)_{t\ge 0}$ on $H^{-1}$ such that $(\ref{ADD})$ holds.  \end{lem}

\beg{proof} It suffices to verify assumptions (K), (H1), (H2), (H3) and (H4) in \cite[Theorem 2.1]{RRW07}.  To verify these assumptions, let $V=L^{1+r}(\mu), H=H^{-1},
V^*$ be the dual space of $V$ w.r.t. $H^{-1}$, and $K= L^{1+r}(\OO\times [0,T]\times E; \P\times\d t\times\mu)$. Then assumption (K) holds for $R(x)= \|x\|_{1+r}^{1+r}, x\in V,$ and $W_1(s)=W_2(s)= s^{1/(1+r)}, s\ge 0.$ Next, from {\bf (A2)} it is easy to see that the hemicontinuity condition (H1) holds for $A(u):= L\Psi(u)+\Phi(u)$, i.e.
$$\R\ni \ll\ni\mapsto \<\Psi(u+\ll v), w\>+ \<\Phi (u+\ll v), (-L)^{-1}w\>,\ \ u,v,w\in V$$ is continuous, where $L$ is understood as its unique extension given by \cite[Lemma 3.3]{RRW07} for $L_{N^*}=L^{\ff{r+1}r}(\mu)$. Moreover, letting $l_0$ be the Lipschitz constant of $\Phi$, we obtain from {\bf (A2)} and {\bf (A3')} that
\beg{equation*}\beg{split}& 2_{V^*}\<A(u)-A(v), u-v\>_V +\|\si_t(u)-\si_t(v)\|_{\L_{HS}(L^2(\mu); H)}^2\\
&\le C\|u-v\|_H^2 -2 \<\Psi(u)-\Psi(v),u-v\>
+ 2 \<\Phi(u)-\Phi(v), (-L)^{-1} (u-v)\>\\
&\le c'\|u-v\|_H^2,\ \ u,v\in V\end{split}\end{equation*}   for some constant $c'>0$. Thus, the weak monotonicity condition (H2) holds.  Again by {\bf(A2)} and {\bf (A3')}, we have \beg{equation*}\beg{split} &2_{V^*}\<A(u),u\>_V +\|\si_t(u)\|_{\L_{HS}(L^2(\mu);H)}^2\\
&\le -2\mu(u\Psi(u))+2 \<\Phi(u), (-L)^{-1}u\> + \ff 1 {\ll_1} \|\si_t(u)\|_{\L_{HS}}^2\\
&\le  c' + c_1 \|u\|_{r+1}^{r+1} +c''\|u\|_H^2\end{split}\end{equation*} for some constants $c', c''>0$. This implies (H3) as $R(u)=\|u\|_{1+r}^{1+r}.$ Finally, {\bf (A2)} also implies that
$$|_{V^*}\<A(v),u\>_V |\le |\mu(\Psi(v)u)|+\|\Phi(v)\|\cdot\|(-L)^{-1}u\|\le c'(1+R(u)+R(v))$$ for some constant $c'>0$ and $R=\|\cdot\|_{1+r}^{1+r}.$ Therefore, (H4) holds.
\end{proof}

Now, let $\E\|X_0\|_{1+r}^{1+r}<\infty.$ To prove that the unique solution $X$ to (\ref{2.1}) is a regular solution in the sense of Definition \ref{D2.1}, we make use of the Galerkin approximations. For any $n\ge 1,$ let $H_n=\text{span}\{e_1,\cdots, e_n\}$. Since $\{e_i\}_{i\ge 1}
\subset L^{1+r}(\mu)$, the orthogonal  projection $\scr P_n: L^2(\mu)\to H_n$ can be  extended to $L^{(1+r)/r}(\mu)$ as
$$\scr P_n u= \sum_{i=1}^n \mu(ue_i) e_i,\ \ \ u\in L^{(1+r)/r}(\mu).$$ Let $\Psi_n(u)=\scr P_n \Psi(u), \Phi_n(u)= \scr P_n \Phi(u)$ for
$u\in L^{1+r}(\mu)$, and let $\si_t^n=\scr P_n \si_t, W_t^n=\scr P_n W_t.$ Finally, let $X_0^n$ be the $L^{1+r}$-best approximation of $X_0$ in $H_n$; that is, $X_0^n$ is the unique $\F_0$-measurable random variable in $H_n$ such that
$$\|X_0-X_0^n\|_{r+1}=\inf_{u\in L^{1+r}(\mu)\cap H_n} \|X_0-u\|_{1+r}.$$ We have $\|X_0^n\|_{1+r}\le 2 \|X_0\|_{1+r}$ and $X_0^{n}\to X_0$ in $L^{1+r}(\mu)$ $\P$-a.s. as $n\to\infty$, see \cite[Theorems 5,6,8]{GS} and \cite[\S 0]{Gess}.

For each $n\ge 1,$ let $X^n= (X_t^n)_{t\in [0,T]}$ be the unique solution to the following finite-dimensional SDE with initial data $X_0^n$:
$$\d X_t^n = \{L\Psi_n (X_t^n)+\Phi_n(X_t^n)\}\d t +\si_t^n(X_t^n)\d W_t^n.$$ Note that since  $H_n$ is an invariant space for $L$, we have $L\Psi_n(u)\in H_n$ for
$u\in L^{1+r}(\mu).$

\beg{lem}\label{L2.2} Assume {\bf (A2)}  and {\bf (A3')} and let $\E\|X_0\|_{1+r}^{1+r}<\infty.$ Then for any $T>0$ there exists a constant $C>0$ independent of $c_1'$ such that
for any $n\ge 1$,
$$\E \sup_{t\in [0,T]}\|X_t^n\|^2 +\sup_{t\in [0,T]}\E\|X_t^n\|_{1+r}^{1+r} + \E \int_0^T \Big\{\EE(\Psi_n(X_t^n),\Psi_n(X_t^n))+c'_1\EE(X_t^n,X_t^n)\Big\}\d t\le C.$$\end{lem}

\beg{proof} Let
$$F(u)= \int_E \d\mu \int_0^u \Psi(s)\d s,\ \ \ u\in H_n.$$ Then $F\in C^2(H_n)$. By the It\^o formula,
we have
\beq\label{T1} \beg{split}&\d F(X_t^n)=\<\Psi(X_t^n), \si_t^n(X_t^n)\d W_t^n\>\\
& + \bigg\{\<\Psi_n(X_t^n), L\Psi_n(X_t^n)+\Phi_n(X_t^n)\>
+\ff 1 2 \sum_{i=1}^n \int_E\Psi'(X_t^n)(\si_t^n (X_t^n)e_i)^2\d\mu\bigg\}\d t.\end{split}\end{equation}
By {\bf (A2)} and {\bf (A3')}, there exists a constant $C_1>1$ independent of $n$ such that for any $u\in H_n$,
\beg{equation*}\beg{split} & \ff 1 {C_1}\|u\|_{1+r}^{1+r}-C_1\le F(u)\le C_1+C_1 \|u\|_{1+r}^{1+r}\\
& |\<\Phi_n(u), \Psi_n(u)\>|\le \|\Phi(u)\|\cdot\|\Psi(u)\|\le C_1+C_1\|u\|_{1+r}^{1+r}\end{split}\end{equation*} and
\beg{equation*}\beg{split} &\E\sum_{i=1}^n \int_E\psi'(X_t^n)(\si_t^n(X_t^n) e_i)^2\d\mu \\
&\le \E\bigg(\|\Psi'(X_t^n)\|_{\ff{1+r}{r-1}}\Big\|\sum_{i=1}^n \Big(
\sum_{k=1}^n \<\si_t(X_t^n) e_i, e_k\>e_k\Big)^2\Big\|_{\ff{1+r}2}\bigg)\le C_1+C_1 \E \|X_t^n\|_{1+r}^{1+r}.\end{split}\end{equation*}
Combining this with (\ref{T1}) we obtain
$$\E\|X_t^n\|_{1+r}^{1+r} \le C_2 +C_2 \int_0^t \E\|X_s^n\|_{1+r}^{1+r}\d s -C_3 \E\int_0^t \EE(\Psi_n(X_s^n),\Psi(X_s^n))\d s,\ \ s\in [0,T]$$ for some constants $C_2, C_3>0$ independent of $n$. This implies
$$\sup_{t\in [0,T]}\E \|X_t^n\|_{1+r}^{1+r} + \E\int_0^T \EE(\Psi_n(X_t^n), \Psi_n(X_t^n))\d t\le C$$ for some constant $C>0$ independent of $n$.

Next, by the It\^o formula,
\beq\label{T2} \d \|X_t^n\|^2 = \big\{2\<X_t^n, L\Psi_n(X_t^n)+\Phi_n(X_t^n)\>+ \|\si_t^n(X_t^n)\|_{\L_{HS}}^2\big\}\d t + 2 \<X_t^n, \si_t^n(X_t^n)\d W_t^n\>.\end{equation}
Since due to {\bf (A2)}   $\Psi'(s)\ge c_1'$, we have
$$\<X_t^n, L\Psi_n(X_t^n)\>=\<X_t^n, L\Psi_n(X_t^n)\> \leq -c_1' \EE(X_t^n, X_t^n).$$ Moreover, by {\bf (A3')} and the Burkholder-Davies inequality for $p=2$ we have
$$\E\sup_{s\in [0,t]} \bigg|\int_0^s \<X_a^n, \si_a^n(X_a^n)\d W_a^n\>\bigg|^2 \le 4\, \E\int_0^t \|\si_s(X_s^n)\|_{\L_{HS}}^2\|X_s^n\|^2\d s\le 4c\,\E\int_0^t\|X_s^n\|^2\d s.$$
Combining this with (\ref{T2}) we conclude that $h_n(t):= \E\sup_{s\in [0,t]} \|X_s^n\|^2$ satisfies
$$ h_n(t)\le C_4+C_4 \int_0^th_n(s)\d s -2c_1' \int_0^t \EE(X_s^n, X_s^n)\d s,\ \ t\in [0,T]$$ for some constant $C_4>0$ independent of $n$.  Therefore,
$$\E\sup_{t\in [0,T]}\|X_t^n\|^2 +c_1'\E\int_0^T \EE(X_s^n,X_s^n)\d s\le C$$   for some constant $C>0$ independent of $n$ and $c_1'$.
\end{proof}

\beg{proof}[Proof of Theorem \ref{T2.1}] To see that the  unique solution $X$ from  Lemma \ref{L2.1} is a regular solution, let us recall the construction of $X$ given in  the proof of \cite[Theorem 2.1]{RRW07}. By Lemma \ref{L2.2} and {\bf (A3')}, there exists a subsequence $n_k\to\infty$, an adapted $ X\in   L^\infty([0,T]\to L^{1+r}(\P\times\mu))$, an adapted
$A\in L^2(\OO\times [0,T]\to H^{-1}; \P\times \d t)$, and some element $Z\in L^2(\OO\times [0,T]\to \L_{HS}; \P\times \d t)$ such that
\beg{enumerate}
\item[(i)] $X^{n_k}\to  X$    $*$-weakly in $L^\infty([0,T]\to L^{1+r}(\P\times\mu)).$
\item[(ii)] $L\Psi_{n_k}(X^{n_k})+\Phi_{n_k}(X^{n_k}) \to A$ weakly in $L^2(\OO\times [0,T]\to H^{-1}; \P\times \d t).$
\item[(iii)] $\si(X^{n_k})\to Z$ weakly in $L^2(\OO\times [0,T]\to \L_{HS}; \P\times \d t)$.
\end{enumerate}
Since these convergence properties  are stronger than those used in the proof of \cite[Theorem 2.1]{RRW07} for $p=2$ and the spaces $K, V, V^*, H$ given in the proof of Lemma \ref{L2.1}, the arguments in the proof of \cite[Theorem 2.1]{RRW07} imply that
 $Z= \si_\cdot(X_\cdot)$ and $A= L\Psi( X)+\Phi( X), \ \P\times\d t$-a.e., and
\beq\label{NN}X_t= \int_0^t A_s\d s +\int_0^t Z_s\d W_s,\ \ t\in [0,T].\end{equation} We are now able to prove the desired regularity properties as follows.

 (a) Since $A\in L^2(\OO\times [0,T]\to H^{-1}; \P\times \d t)$ and $A= L\Psi(X)+\Phi(X)\ \P\times\d t$-a.e., we have
 $L\Psi( X)+\Phi( X)\in L^2(\OO\times [0,T]\to H^{-1}; \P\times \d t).$ Moreover, since $\Phi$ is Lipschitz continuous and $$\|\cdot\|_{1+r}\ge \|\cdot\|\ge \ff 1 {\ss\ll_1}\|\cdot\|_{H^{-1}},$$   Lemma \ref{L2.1} implies  that $\Phi(X)\in  L^2(\OO\times [0,T]\to H^{-1}; \P\times \d t).$ Therefore, $L\Psi( X)\in L^2(\OO\times [0,T]\to H^{-1}; \P\times \d t),$ that is,
\beq\label{ABC}\E\int_0^T \EE(\Psi(X_t),\Psi(X_t))\d t = \E\int_0^T \|L\Psi(X_t)\|_{H^{-1}}^2 <\infty.\end{equation}Since
   {\bf (A2)}   implies  $\Psi'\ge c_1 $ so that $\EE(X_s,\Psi(X_s))\ge c_1\EE(X_s,X_s)$, it follows from (\ref{ABC}) that
  $$c_1'\E\int_0^T \EE(X_t,X_t)\d t<\infty.$$

 (b) When $\Psi(s)=c |s|^{r-1}s$ and $ \Phi(s) =c' s$  for $c>0,$ the right continuity of the solution in $L^2(\m)$ and
 $\E \sup_{t\in [0,T]} \|X_t\|^2<\infty$ are  ensured by \cite[Theorem 1.2(4)]{RW08}. Let $c_1'>0$, so that  $X\in L^2(\OO\times [0,T]\to H^1; \P\times\d t)$. To see that $X_t$ is continuous in $L^2(\mu)$, we make use of \cite[Theorem A.2]{RRW07}.
  Let now $K= L^2(\OO\times [0,T]\to H^1; \P\times \d t), H=L^2(\mu), V= H^1$ and $V^*= H^{-1}.$   Then the condition (K) in \cite{RRW07} holds for $R(u)=\EE(u,u)
  =\|u\|_{H^1}^2$ and $W_1(s)=W_2(s)=\ss s, s\ge 0.$ Since $A\in K^*:= L^2(\OO\times [0,T]\to H^{-1}; \P\times \d t)$ and $Z\in J:= L^2(\OO\times [0,T]\to \scr L_{HS}; \P\times \d t)$ (see (2.8) in \cite{RRW07}),   according to \cite[Theorem A.2]{RRW07}, (\ref{NN}) implies that $X_t$ is continuous in $H(=L^2(\mu))$ such that $\E\sup_{t\in [0,T]}\|X_t\|^2<\infty$ and that the It\^o formula
  $$\|X_t\|^2=\|X_0\|^2+\int_0^t\Big\{2_{H^1}\<X_s, A_s\>_{H^{-1}}+\|\si_s(X_s)\|_{\L_{HS}}^2\Big\}\d s  +2 \int_0^t \<\si_s(X_s)\d W_s, X_s\>,\ \ t\in [0,T],$$ holds. This coincides with the desired It\^o formula since $A= L\Psi(X)+\Phi(X)\ \P\times\d t$-a.e.

 (c)  It remains to show that $\sup_{t\in [0,T]} \E\|X_t\|_{1+r}^{1+r}<\infty$ for cases (1) and (2). Since $X\in L^\infty([0,T]\to L^{1+r}(\P\times\mu)),$ there exists a constant $C>0$ such that $\E\|X_t\|_{1+r}^{1+r}\le C$ holds $\d t$-a.e. Since $X_t$ is right-continuous in $L^2(\mu)$, this and the Fatou lemma imply that for any $t\in [0,T]$,
 \beg{equation*}\beg{split} \E \|X_t\|_{r+1}^{r+1} &=  \E \sup_{u} |\mu(X_tu)|^{1+r} =  \E\sup_u \liminf_{s\downarrow t}|\mu(X_su)|^{1+r}\\
 &\le \liminf_{s\downarrow t} \E \sup_u |\mu(X_su)|^{1+r}\le   \liminf_{s\downarrow t} \E \|X_s\|_{1+r}^{1+r} \le C,\end{split}\end{equation*}
 where $\sup$ is taken over all $u\in L^2(\mu)$  with $ \|u\|_{\ff r{1+r}}\le 1$.
 \end{proof}

\subsection{Comparison theorem and $L^1$-Lipschitz continuity}

In this subsection we consider the following equation with additive noise:
\beq\label{1.1'} \d X_t= \big\{L\Psi( X_t)  +  \Phi(X_t)\big\}\d t +\si_t\d W_t,\end{equation} where $\si_t, \Psi$ and $\Phi$ satisfy
{\bf (A2)} and {\bf (A3}). Let $\tt\Phi$ be another Lipshitz continuous function. We shall compare regular solutions to (\ref{1.1'}) with those to the equation
\beq\label{PP} \d\tt X_t= \big\{L\Psi(\tt X_t)  +\tt \Phi(\tt X_t)\big\}\d t +\si_t\d W_t.\end{equation}

\beg{thm}\label{T2.2} Assume {\bf (A2)}, {\bf (A3)} and let $\tt\Phi\le \Phi.$ Let $X_t$ and $\tt X_t$ be solutions in the sense of Definition $\ref{D2.1}$ to $(\ref{1.1'})$ and $(\ref{PP})$ respectively. If either $c_1'>0,$ or $\Psi(s)= c |s|^{r-1}s, \Phi(s)=c's$ for some $c>0$ and $c'\in\R$, then these solutions are regular and   $\P$-a.s. $\tt X_0\le X_0$ implies $\P$-a.s. $\tt X_t\le X_t$ for all $t\in [0,T]. $ \end{thm}

Let us first explain the main idea of the proof. The regularity of the solutions follows from Theorem \ref{T2.1}. To prove $\tt X_t\le X_t$, let $h_k\in C^1_b(\R)$ such that $h_k'\ge 0,$ $0\le h_k\le 1, h_k(s)=0$ for $s\le 0$, and $h_k\to 1_{(0,\infty)}$ as $k\to\infty.$ By the definition of regular solutions, $\P$-a.s.,
$\tt X_t- X_t$ is $\d t$-a.e. differentiable in $H^{-1}$ with
$$\ff{\d }{\d t} (\tt X_t-X_t)= L\big\{\Psi(\tt X_t)-\Psi(X_t)\big\} +\tt\Phi(\tt X_t)-\Phi(X_t).$$ Moreover,
$h_k\big(\Psi(\tt X_t)-\Psi(X_t)\big)\in H^1\ \P\times\d t$-a.e. Therefore, noting that $\Phi$ is Lipschitzian and $\tt\Phi\le \Phi$, we have $\P\times\d t$-a.e.,
\beg{equation}\label{ZZ}\beg{split} &\int_E \Big\{h_k\big(\Psi(\tt X_t)-\Psi(X_t)\big)\ff{\d }{\d t} (\tt X_t-X_t)\Big\}\d\mu\\
&=-\int_E\EE\big(h_k\big(\Psi(\tt X_t)-\Psi(X_t)\big), \Psi(\tt X_t)-\Psi(X_t)\big)\d\mu\\
&\qquad + \int_E h_k\big(\Psi(\tt X_t)-\Psi(X_t)\big)\cdot \big(\tt\Phi(\tt X_t)-\Phi(X_t)\big) \d \mu\\
&\le l_0\int_E h_k\big(\Psi(\tt X_t)-\Psi(X_t)\big)\cdot |\tt X_t-X_t|\d\mu,\end{split}\end{equation} where $l_0$ is the Lipschitz constant of $\Phi$.
By letting $k\to\infty$ we may write formally
\beg{equation}\label{ZZD}``\ \ff{\d }{\d t} \mu\big((\tt X_t-X_t)^+\big) = \int_E \Big\{1_{\{\tt X_t>X_t\}}\ff{\d }{\d t} (\tt X_t-X_t)\Big\}\d\mu\ \text{''}
 \le l_0 \mu((\tt X_t-X_t)^+),\end{equation}
and hence, $(\tt X_t-X_t)^+=0$ if $\tt X_0\le X_0$ as desired. The last step is however not rigorous since $\ff{\d}{\d t}(\tt X_t-X_t)$ exists only in $H^{-1}$ so that the terms in $``$\ ... '' do not make sense in general.  To make the argument rigorous, we consider the following approximating equations for $\vv\in (0,1)$:
\beg{equation}\label{XX}\beg{split} &\d X_t^\vv = \big\{(1-\vv L)^{-1}L \Psi(X_t^\vv)+\Phi(X_t^\vv)\big\}\d t +\si_t\d W_t,\ \ X_0^\vv=X_0,\\
&\d \tt X_t^\vv = \big\{(1-\vv L)^{-1}L \Psi(\tt X_t^\vv)+\Phi(\tt X_t^\vv)\big\}\d t +\si_t\d W_t,\ \ \tt X_0^\vv=\tt X_0.\end{split}\end{equation}

\beg{lem}\label{L2.3} We have
$$\lim_{\vv\to 0} \E\int_0^T \Big(\|X_t-X_t^\vv\|_{1+r}^{1+r}+\|\tt X_t-\tt X_t^\vv\|_{1+r}^{1+r}\Big)\d t=0.$$\end{lem}
\beg{proof} We only consider the limit for $X_t-X_t^\vv$.  Since $\Phi$ is Lipschitz continuous, there exists a constant $C_1>0$ independent of $\vv$ such that
$$\<\Phi(X_t)-\Phi(X_t^\vv), X_t-X_t^\vv\>_{H^{-1}}\le C_1 \|X_t-X_t^\vv\|\cdot\|X_t-X_t^\vv\|_{H^{-1}}.$$ Moreover, from the proof of Lemma \ref{L2.1} we see that
$$\sup_{\vv\in (0,1)} \E\int_0^T (\|X_t\|_{1+r}^{1+r}+\|X_t^\vv\|_{1+r}^{1+r})\d t <\infty.$$ Combining this with the growth condition   $|\Psi(s)|\le c'(1+|s|^r)$ ensured by {\bf (A2)}, we obtain
\beg{equation*}\beg{split} &\E\int_0^T  \big| \<(1-(1-\vv L)^{-1})(X_t-X_t^\vv), \Psi(X_t^\vv)\>\big| \d t \\
&\le \E\int_0^T  \|\Psi(X_t^\vv)\| \cdot\|X_t-X_t^\vv\|\d t\le C\end{split}\end{equation*} for
a constant $C>0$ independent of $\vv\in (0,1).$
Therefore, by {\bf(A2)} and the It\^o formula,
\beg{equation*}\beg{split} &\E \|X_t-X_t^\vv\|_{H^{-1}}^2 \le -2 c_1 \int_0^t \E   \|X_s-X_s^\vv\|_{1+r}^{1+r}   \d s
  +C_1 \int_0^t \E   \|X_s-X_s^\vv\|_{H^{-1}}^2\d s \\
  &\qquad  + 2\vv \int_0^t \E \Big(|\<(1-(1-\vv L)^{-1})(X_s-X_s^\vv), \Psi(X_s^\vv)\>|+|\<(1-\vv L)^{-1}(X_s-X_s^\vv),X_s^\vv\>|\Big)\d s\\
&\le C_2\int_0^t \E \|X_s-X_s^\vv\|_{H^{-1}}^2\d s -2 c_1 \int_0^t \E \|X_s-X_s^\vv\|_{1+r}^{1+r}\d s + 2C_2\vv,\ \ t\in [0,T] \end{split}\end{equation*}  holds for some constants $C_1,C_2>0.$  This implies $\lim_{\vv\to 0}\E\int_0^T \|X_t-X_t^\vv\|_{1+r}^{1+r}\d t=0.$
\end{proof}

\beg{proof}[Proof of Theorem \ref{T2.2}] Since $(1-\vv L)^{-1}L$ is a bounded operator for any $\vv>0$, the associated Dirichlet space and its dual space w.r.t. $L^2(\mu)$ coincide with $L^2(\mu).$ So, by Definition \ref{D2.1},
$\ff{\d (\tt{X_t}^\vv -X_t^\vv)}{\d t}$ exists in $L^2(\OO\times [0,T]\to L^2(\m); \P\times\d t)$. We aim to prove, instead of (\ref{ZZD}), that
\beq\label{ZZD'}\ff{\d }{\d t} \mu\big((\tt X_t^\vv-X_t^\vv)^+\big) = \int_E \Big\{1_{\{\tt X_t>X_t\}}\ff{\d }{\d t} (\tt X_t^\vv-X_t^\vv)\Big\}\d\mu\le l_0 \mu(\tt X_t^\vv-X_t^\vv)^+),\end{equation} which implies
$\tt X_t^\vv\le X_t^\vv$ for all $t\in [0,T] $ since $\tt X_0^\vv\le X_0^\vv.$ Firstly,  replacing $(\Psi,X_t,\tt X_t)$ in  (\ref{ZZ}) by $(\Psi,X_t^\vv, \tt X_t^\vv)$  and letting $k\to\infty$, we obtain the inequality in (\ref{ZZD'}). To verify the   equality in (\ref{ZZD'}), we note that
 $$\ff{\d }{\d t} (\tt X_t^\vv-X_t^\vv)=(1-\vv L)^{-1}L\big(\Psi(\tt X_t^\vv)-\Psi(X_t^\vv)\big)+\tt\Phi(\tt X_t)-\Phi(X_t^\vv) $$ and   $(1-\vv L)^{-1}L=\ff 1 \vv (1-\vv L)^{-1}-\ff 1 \vv$ on $L^1(\mu)$ imply
 \beg{equation*}\beg{split} &\sup_{0\le s <t\le T}\bigg|\ff{(\tt X_t^\vv-X_t^\vv)^+-(\tt X_s^\vv-X_s)^+}{t-s}\bigg|\\
 &\le \ff 1 \vv\big(1+(1-\vv L)^{-1}) \sup_{r\in [0,T]}\Big\{|\Psi(\tt X_r^\vv)|+ |\Psi(X_r^\vv)|+ |\tt\Phi(\tt X_r^\vv)|+|\Phi (X_r^\vv)|\Big\}. \end{split}\end{equation*}  By   the contraction property of  $(1-\vv L)^{-1}$    on $L^1(\mu)$, Lemma \ref{L2.1}, and the growth conditions on $\Psi,\Phi,\tt\Phi$, we see that the upper bound  is in $L^1(\mu)$. Therefore, the equality in (\ref{ZZD'}) follows from the dominated convergence theorem with $s\to t$.

 Now, by $\tt X_t^\vv\le X^\vv_t$, we have
 $$\E\int_0^T \mu\big((\tt X_t^\vv-X_t^\vv)^+\big)\d t=0,\ \ \vv\in (0,1).$$ Letting $\vv\to 0$ and using Lemma \ref{L2.3}, we arrive at
 $$\E\int_0^T \mu\big((\tt X_t-X_t)^+\big)\d t=0.$$ Therefore, $\tt X_t\le X_t$ holds $\P\times\d t\times\mu.$ Since  due to Theorem \ref{T2.1} $X_t$ and $\tt X_t$ are
 right-continuous in $L^2(\mu)$, we conclude that $\P$-a.s., $\tt X_t\le X_t$ in $L^2(\mu)$ holds for all $t\in [0,T].$
\end{proof}

Next, we have the following $L^1$-Lipschitz continuity w.r.t. initial data of the solutions.

\beg{thm}\label{T2.3} Assume {\bf (A2)},{\bf (A3)} and $(\ref{PH})$. We have
 $$\|X_t(x)-X_t(y)\|_1\leq \e^{Kt}\|x-y\|_1, \quad x, y\in L^{1+r}(\mu).$$
 \end{thm}
\beg{proof}  Let $X_t^\vv(x)$ be as in (\ref{XX}) for $X_0=x\in L^{1+r}(\mu)$. Repeating the proof of Theorem \ref{T2.2} with $(\tt X^\vv, X^\vv)$ replaced by
$(X^\vv (x), X^\vv (y))$, we obtain
$$\d\|(X_t^\vv(x)-X_t^\vv(y))^+\|_1=\<1_{\{X_t^\vv(x)-X_t^\vv(y)>0\}}, \d (X_t^\vv(x)-X_t^\vv(y))\>\d t \le K \|(X_t^\vv(x)-X_t^\vv(y))^+\|_1\d t.$$ Then
$$\|(X_t^\vv(x)-X_t^\vv(y))^+\|_1\le \e^{Kt}\|x-y\|_1.$$ The same holds by switching $x$ and $y$ so that
$$\|X_t^\vv(x)-X_t^\vv(y)\|_1\le \e^{Kt}\|x-y\|_1,\ \ t\ge 0.$$  Since due to Lemma \ref{L2.3} there exists a sequence $\vv_n\downarrow 0$ such that for any $T>0$
$$\lim_{n\to\infty}\int_0^T\big(\|X_t^{\vv_n}(x)-X_t(x)\|_{1+r}^{1+r}+ \|X_t^{\vv_n}(y)-X_t(y)\|_{1+r}^{1+r}\big)\d t=0,$$ this implies that
$\|X_t(x)-X_t(y)\|_1\le \e^{Kt}\|x-y\|_1$ holds $\d t$-a.e. Then the proof is finished by the continuity of the solutions. \end{proof}

\subsection{Riesz-Markov representation theorem}
Let $\tt E$ be a locally compact separable metric space so that $\si(C_0(\tt E))=\tt {\scr B}$ (the Borel $\si$-field on $\tt E$), and let $\tt\C\subset C_0(\tt E)$ be a subspace  such that the following assumption holds:

\paragraph{(A)} \ for any $f\in C_0(\tt E)$, there exists $\tt f\in \tt\C$ such that for any $\vv>0$, there exists $f_\vv\in \tt\C$

\ \ \ \ \ such that
$|f-f_\vv|\le\vv \tt f.$

\ \newline Let $C_0^+(\tt E)$ and $\tt\C^+$ denote the classes of non-negative elements in $C_0(\tt E)$ and $\tt\C$ respectively.
\beg{thm}\label{T2.3} Assume {\bf (A)}. For any positive linear functional
$\LL: \tt\C\to \R$, there exists a unique measure $\mu$ on $\tt E$ such that
\beq\label{LL} \mu(f):=\int_{\tt E}f\d\mu = \LL(f),\ \ f\in \tt\C.\end{equation}
\end{thm}

\beg{proof} (a) The uniqueness. Let $\mu$ and $\tt\mu$ be two measures satisfying (\ref{LL}), then for any $f\in C_0(\tt E)$, and for  $\tt f$ and $f_\vv$ in {\bf (A)}, we have $f_\vv+\vv f\in \tt\C$ so that
$$ \mu(f)\le \mu(f_\vv+\vv\tt f)= \tt\mu (f_\vv+\vv \tt f) \le  \tt\mu(f)+ 2\vv \LL(\tt f).$$ Letting $\vv\to 0$ we obtain $\mu(f)\le\tt\mu(f)$. Similarly, $\tt\mu(f)\le \mu(f).$ Therefore, $\mu=\tt\mu$.

(b) The existence. For any $f\in C_0^+(\tt E)$, let
$$\bar\LL(f)=\sup\{\LL(g):  g\le f, g\in \tt\C\}.$$ Since $0\in \tt\C$, we have $\bar\LL(f)\ge 0$ for $f\in C_0^+(\tt E).$ Next, it is easy to see that  $\bar\LL$ is increasing monotone and $\bar\LL=\LL$ holds on $\tt\C^+.$
Moreover, by {\bf (A)}, for $f\in C_0^+(\tt E)$ there exists $\tt f,g\in \tt\C$ such that $|f-g|\le \tt f.$ Then $\tt f+g\in \tt \C^+$ so that
$$ \LL(f)\le \bar\LL(\tt f+g) =\LL(\tt f+g)<\infty.$$ Therefore, letting $\bar\LL(f)= \bar\LL(f^+)-\bar\LL(f^-)$, we extend $\bar\LL$ to a finite positive functional on $C_0(\tt E)$ such that $\bar\LL=\LL$ holds on $\tt\C$. Then it suffices to show that
\beq\label{2C2}  \LL(f+g)= \LL(f)+\bar\LL(g),\ \ f,g\in C_0(\tt E).\end{equation} Indeed, it is trivial to see that
$\bar\LL(cf)=c\bar\LL(f)$ for $f\in C_0(\tt E)$ and $c\in\R$. Then (\ref{2C2}) implies that $ \LL: C_0(\tt E)\to \R$
is a positive linear functional. By the Riesz-Markov representation theorem,  there exists a  unique locally finite  measure $\mu$ on $\tt E$ such that
$$\mu(f)= \bar\LL(f),\ \ f\in C_0(\tt E).$$ Since $\LL(f)=\bar\LL(f)$ holds for $f\in \tt\C$, this implies (\ref{LL}).

Now, let $f,g\in C_0(\tt E).$ By {\bf (A)}, there exist $\tt f,\tt g\in \tt\C$ such that for any $\vv>0$ there exist $f_\vv, g_\vv\in \tt\C$ such that $|f-f_\vv|\le\vv \tt f, |g-g_\vv|\le \vv \tt g.$ We have
\beg{equation*}\beg{split} &\bar\LL(f+g)\le \bar\LL(f_\vv+g_\vv+\vv \tt f+\vv\tt g)\\
&= \LL(f_\vv-\vv \tt f)+\LL(g_\vv-\vv \tt g)+ 2\vv\LL(\tt f+\tt g)
\le\bar\LL(f)+\bar\LL(g)+2\vv\LL(\tt f+\tt g),\end{split}\end{equation*} and conversely,
\beg{equation*}\beg{split}&\bar\LL(f)+\bar\LL(g)\le \LL(f_\vv+\vv \tt f)+\LL(g_\vv+\vv \tt g)\\
& = \LL(f_\vv+g_\vv -\vv\tt f-\vv \tt g)+2\vv\LL(\tt f+\tt g)\le
\bar\LL(f+g)+2\vv \LL(\tt f+\tt g).\end{split}\end{equation*} Letting $\vv\to 0$ we prove (\ref{2C2}).
\end{proof}

\section{Proof of Theorem \ref{T1.1}}

For $n\ge 1$, let $$\Phi^{(n)} (s)=\Phi(s)+n s^-,\ \ \ s\in\R.$$  Consider the following penalized  equation:
\beq\label{0.1}
  \d X_t^{(n)} = L\Psi(X_t^{(n)})dt+ \Phi^{(n)}(X_t^{(n)})\d t+\si_t \d W_t,\ \ X_0^{(n)}=X_0.
\end{equation} Let $$\nu^{(n)}_t(\d z)= n \bigg(\int_0^t(X_s^{(n)}(z))^{-}\d s\bigg)\mu(\d z),\ \ n\ge 1.$$
By Theorem \ref{T2.1}, for each $n$, this equation has a unique regular solution in the sense of Definition \ref{D2.1}.
 We will show that    $( X_t, \nu_t)=\lim_{n\to\infty}(X_t^{(n)}, \nu^{(n)})$ exists and  gives rise to a solution to equation (\ref{1.1}) via (\ref{D1}).

\subsection{Construction and properties of  $ X$}

 By Theorem \ref{T2.1} and Theorem \ref{T2.2}, $\{X^{(n)}\}_{n\ge 1}$ is an increasing sequence of continuous adapted processes in $L^2(\mu)$ such that
 \beq\label{3.1} X^{(n)}, \Psi(X^{(n)})\in L^2(\OO\times [0,T]\to H^1; \P\times \d t),\ \ X^{(n+1)}\ge X^{(n)},\ \ n\ge 1.\end{equation}
 Let $$ X=\lim_{n\to\infty}X^{(n)}.$$

 \beg{lem}\label{L3.1}   $X^{(n)}\to  X$ in $L^2(\OO\times [0,T]\times E; \P\times \d t\times\mu)$ and
 \beq\label{*1}\E\sup_{t\in [0,T]}\| X_t\|^2 <\infty.\end{equation} Consequently, $X_t^{(n)}\to  X_t$ holds in $L^2(\OO\times E;\P\times \mu)$ for all $t\in [0,T].$\end{lem}

 \beg{proof}  By the It\^o formula in Theorem \ref{T2.1} and using {\bf (A2)}, {\bf (A3)}, we obtain
 \beg{equation*}\beg{split} \d\|X_t^{(n)}\|^2 &\leq \big\{2\<\Phi(X_t^{(n)}, X_t^{(n)}\> + 2n \<(X_t^{(n)})^-, X_t^{(n)}\> \\ &\qquad \qquad  +\|\si_t\|_{\L_{LS}}^2\big\}\d t + 2 \<\si_t\d W_t, X_t^{(n)}\>\\
 &\le \big\{C_1 +C_1 \|X_t^{(n)}\|^2  \big\}\d t + 2 \<\si_t\d W_t, X_t^{(n)}\>\end{split}\end{equation*} for some constant
 $C_1>0$ independent of $n$. As shown in the second part in the proof of Lemma \ref{L2.2}, this implies
 \beq\label{3.2} \E\sup_{t\in [0,T]}\|X_t^{(n)}\|^2 \le C\end{equation} for some constant independent of $n$. Noting that
 $$\|(X_t^{(n)})^+\|\uparrow \| X_t^+\|,\ \ \| X_t^-\|\le \|(X_t^{(1)})^-\|,$$ this implies
\beg{equation*}\beg{split} \E\sup_{t\in [0,T]}\| X_t\|^2 &\le \E \sup_{t\in [0,T]}\sup_{n\ge 1}\|(X_t^{(n)})^+\| +\E\sup_{t\in [0,T]}\|X_t^{(1)}\|^2\\
& = \lim_{n\to\infty} \E \sup_{t\in [0,T]}\|(X_t^{(n)})^+\| +\E\sup_{t\in [0,T]}\|X_t^{(1)}\|^2\le 2C.\end{split}\end{equation*}
 Since $X^{(n)}\uparrow  X,\ \P\times\d t\times\mu$-a.e. and $|X^{(n)}|\le  X^+ +(X^{(1)})^-$, by the dominated convergence theorem we conclude that
 $X^{(n)}\to  X$ in $L^2(\OO\times [0,T]\times E;\P\times\d t\times \mu).$
 \end{proof}

 \beg{lem}\label{L3.2} $\Psi(X^{(n)})\to \Psi( X)$ weakly in $L^2(\OO\times [0,T]\to H^1;\P\times\d t)$ and
 $$\sup_{t\in [0,T]}\E\| X_t\|_{1+r}^{1+r}+\E\int_0^T \EE(\Psi( X_t),\Psi( X_t))\d t<\infty.$$ Moreover, $ X\ge 0,\ \P\times\d t\time\mu$-a.e.
 \end{lem}

 \beg{proof} For $m\ge 1$, let $\phi_m\in C_b^\infty(\R)$ such that $0\le \phi_m'\le 2$ and
 $$\phi_m(s)=\beg{cases} s,\ &\text{if}\ |s|\le m,\\
 m+1,\ &\text{if}\ s\ge m+1,\\
 -m-1, \ &\text{if}\ s\le -m-1.\end{cases}$$ Define
 $$F_m(u)=\int_E\d\mu \int_0^u \Psi\circ\phi_m(s)\d s,\ \ \ u\in L^2(\mu).$$ Then $F_m\in C_b^2(L^2(\mu))$ with
 \beg{equation*}\beg{split}&\pp_{v_1} F_m(u)= \int_E \Psi(\phi_m(u)) v_1\d\mu,\\
 &\pp_{v_1}\pp_{v_2}F(u)= \int_E \Psi'\circ\phi_m(u) \phi_m'(u) v_1v_2\d\mu,\ \ u,v_1,v_2\in L^2(\mu).\end{split}\end{equation*} Since due to Theorem \ref{T2.1} we have $X^{(n)}, \Psi(X^{(n)})\in L^2(\OO\times[0,T]\to H^1),$ by the It\^o formula we obtain
 \beq\label{3A} \beg{split} &\d F_m (X_t^{(n)} ) = \bigg\{\<\Phi(X_t^{(n)}), \Psi\circ\phi_m (X_t^{(n)})\> -\EE\big(\Psi(X_t^{(n)}), \Psi\circ \phi_m(X_t^{(n)})\big)\\ & +n \<(X_t^{(n)})^-, \Psi\circ\phi_m(X_t^{(n)})\>
+\ff 1 2 \sum_{i=1}^\infty \int_E (\si_te_i)^2 \Psi'\circ\phi_m(X_t^{(n)})
 \phi_m'(X_t^{(n)})\d\mu\bigg\} \d t\\
 &+ \<\si_t\d W_t, \Psi\circ\phi_m(X_t^{(n)})\>.\end{split}\end{equation} Since $\Psi(0)=0, \phi_m(s)\ge 0$ for $s\le 0$, and $\Psi'\ge 0$ imply $s^-\Psi\circ\phi_m(s)\le 0$, we have $$\<(X_t^{(n)})^-, \Psi\circ\phi_m(X_t^{(n)})\>\le 0.$$ Combining this with (\ref{3A}) and the property of $\phi_m$, we obtain
\beg{equation*}\beg{split}  &F_m(X_{t_2}^{(n)})-F_m(X_{t_1}^{(n)}) \le \int_{t_1}^{t_2} \Big\{\<\Phi(X_t^{(n)}), \Psi\circ\phi_m (X_t^{(n)})\> -\EE\big(\Psi(X_t^{(n)}), \Psi\circ \phi_m(X_t^{(n)})\big)\Big\}\d t\\
&+C_1 \sum_{i=1}^\infty \int_{t_1}^{t_2}\mu\big((\si_te_i)^2(1+|X_t^{(n)}|^{r-1})\big)\d t + \int_{t_1}^{t_2}\<\si_t\d W_t, \Psi\circ \phi_m(X_t^{(n)})\>\end{split}\end{equation*} for some constant $C_1>0$ independent of $m,n$, and all $ 0\le t_1\le t_2\le T$. Letting $m\to\infty$
we arrive at
\beg{equation*}\beg{split} \d F(X_t^{(n)}) \le  &\Big\{\<\<\Phi(X_t^{(n)}), \Psi (X_t^{(n)})\> -\EE\big(\Psi(X_t^{(n)}), \Psi(X_t^{(n)})\big)\\
 &\quad +C_1 \sum_{i=1}^\infty \int_{t_1}^{t_2}\mu\big((\si_te_i)^2(1+|X_t^{(n)}|^{r-1})\Big\}\d t
+ \<\si_t\d W_t, \Psi(X_t^{(n)})\>,\end{split}\end{equation*} where $F$ is defined as in the proof of Lemma \ref{L2.2}. Therefore, by repeating the proof
 of Lemma \ref{L2.2}, we obtain
 \beq\label{**} \sup_{t\in [0,T]}\E\|X_t^{(n)}\|_{1+r}^{1+r}+ \E\int_0^T \EE(\Psi(X_t^{(n)}), \Psi(X_t^{(n)}))\d t\le C\end{equation} for some constant $C>0$ independent of $n$. Since $X^{(n)}\uparrow  X$ and $\Psi$ is increasing and continuous, we have $\Psi(X^{(n)})\uparrow \Psi( X)$. Therefore, as in the proof of Lemma \ref{L3.1} we see that (\ref{**}) implies that $\Psi(X^{(n)})\to \Psi( X)$ in $L^2(\OO\times [0,T]\times \E;\P\times\d t\times \mu)$ and (at least a subsequence thereof) weakly in $L^2(\OO\times [0,T]\to H^1)$, as well as
 $$\sup_{t\in [0,T]}\E\| X_t\|_{1+r}^{1+r}+ \E\int_0^T \EE(\Psi( X_t),\Psi( X_t))\d t<\infty.$$

 Finally, we prove that $ X\ge 0$. To this end, let $u\in H^1$. Since $\Psi(X^{(n)})\to \Psi( X)$ weakly in $L^2(\OO\times [0,T]\to H^1)$,  $X^{(n)}\to X$ in $L^2(\OO\times [0,T]\times \E;\P\times\d t\times \mu)$, and $\Phi$ is Lipschitz continuous, we have, in $L^1(\P)$
\beg{equation*}\beg{split} &\int_0^T \< X_t^-,u\>\d t  =\lim_{n\to\infty} \int_{[0,T]\times E} u(X_t^{(n)})^-\d t\d\mu \\
&=\lim_{n\to\infty}\ff 1 n \bigg\{ \<u,  X_T-X_0\>  -\int_0^T \Big\{\<\Phi(X_t^{(n)}), u\>-\EE(\Psi(X_t^{(n)}), u)\Big\}\d t
-\int_0^T\<\si_t\d W_t, u\>\bigg\}\\
&=0.\end{split}\end{equation*} Since $H^1$ is dense in $L^2(\mu)$, this implies that $\int_0^T  X_t^-\d t=0$ in $\P\times\mu$-a.e. Therefore,
$ X\ge 0, \P\times \d t\times \mu$-a.e.
 \end{proof}

\subsection{Construction and properties of  $\eta$}

\beg{lem}\label{L3.3} As $n\to\infty,$ $\nu_t^{(n)}$ converges vaguely to some locally bounded random measure $\bar\nu_t$ on $E$ such that $\P$-a.s.
\beq\label{S'} \beg{split}\bar\nu_t(f)= &\<f,  X_t-X_0\>-\bigg\<f, \int_0^t\si_s\d W_s\bigg\>\\
&+\int_0^t\Big\{\EE(f, \Psi( X_s))
 -\<f, \Phi( X_s)\>\Big\}\d s ,\ \ f\in H^1\cap C_0(E).\end{split}\end{equation}Consequently, $\bar\nu_t$ is an adapted increasing process on $\scr M_c$. \end{lem}

\beg{proof}
 By (\ref{0.1}) and noting that $\{\Psi(X^{(n)})\}_{n\ge 1}$ is a bounded sequence in $L^2(\OO\times [0,T]\to H^1; \P\times \d t)$,    we have $\P$-a.s.
 \beg{equation*}\beg{split} &\nu^{(n)}_t(f)= n\int_E f(z)\mu(\d z)\int_0^t(X^{(n)}_s(z))^-\d s  \\
 &= \<f, X_t^{(n)}- X_0\> -\bigg\<f, \int_0^t\si_s\d W_s\bigg\>
  +\int_0^t\Big\{\EE(f, \Psi(X_s^{(n)}))
 -\<f, \Phi(X_s^{(n)})\>\Big\}\d s \end{split}\end{equation*} for all   $f\in H^1, t\in [0,T]$.  According to Lemmas \ref{L3.1} and \ref{L3.2}, selecting a subsequence if necessary,  we conclude that $\P$-a.s.
 \beq\label{*B} \beg{split}&\LL_t(f):= \lim_{n\to\infty} \nu^{(n)}_t(f)\\
 &=  \<f,  X_t-X_0\>-\bigg\<f, \int_0^t\si_s\d W_s\bigg\>
+\int_0^t\Big\{\EE(f, \Psi( X_s))
 -\<f, \Phi( X_s)\>\Big\}\d s \end{split}\end{equation} exists for all $f\in H^1$. Since $\nu_t^{(n)}\ge 0$,  this implies that
 $\LL_t: H^1\cap C_0(E)\to\R$ is a (random) positive linear functional.  By {\bf (A1)} and Theorem \ref{T2.3}, there exists a unique locally bounded (random) measure $\bar\nu_t$ on $E$ such that
 $$\bar\nu_t(f):=\int_E f(z)\bar\nu_t(\d z)= \LL_t(f),\ \ f\in H^1\cap C_0(E).$$

 Next, to see that $\nu_t^{(n)}\to \bar\nu_t$ vaguely, we first note that (\ref{*B}) and (\ref{S'}) imply
 \beq\label{S''} \lim_{n\to\infty} \nu_t^{(n)}(f)= \bar\nu_t(f),\ \ \ f\in H^1\cap C_0(E).\end{equation}
 Now, let $f\in C_0(E)$. By {\bf (A1)}, there exists $\tt f \in H^1\cap C_0(E)$ such that for any $\vv>0$, $|f- f_\vv|\le \vv\tt f $ holds for some  $f_\vv\in H^1\cap C_0(E)$. Then
 $$\limsup_{n\to\infty} |\nu_t^{(n)}(f)-\bar\nu_t(f)|\le \limsup_{n\to\infty}\big\{|\nu_t^{(n)}(f_\vv)-\bar\nu_t(f_\vv)|+\vv (\nu_t^{(n)}(\tt f)+\bar\nu_t(\tt f))\big\}= 2\vv \bar\nu_t(\tt f).$$ Letting $\vv\to 0$ we conclude that $\lim_{n\to\infty}\nu_t^{(n)}(f)=\bar\nu_t(f).$ Since $\{\bar\nu_t^{(n)}\}_{n\ge 1}$ are locally finite measures, $\bar\nu_t: C_0(E)\to \R$ is a non-negative linear functional and thus is realized by a locally finite measure according to the Riesz-Markov representation theorem, denoted again by $\bar\nu_t$.
 Finally, since $\nu_t^{(n)}$ is increasing in $t$, so is $\bar\nu_t$.
\end{proof}

To construct $\eta$, we observe from {\bf (A3)}, the Lipschitz continuity of $\Phi$ and Lemma \ref{L2.3} that, (\ref{S'}) provides a
bounded linear functional $\bar \nu_t: C_0(E)\cap H^1\to \R$. Since the Dirichlet form is regular, $C_0(E)\cap H^1$ is dense in $H^1$, it can $\P$-a.s. be uniquely extended to an element $\eta_t\in H^{-1}$ such that $\P$-a.s.
\beq\label{DDE} \beg{split} _{H^{-1}}\<\eta_t,f\>_{H^1}=& \<f, X_t-X_0\>-\bigg\<f, \int_0^t\si_s\d W_s\bigg\>\\
&+\int_0^t\Big\{\EE(f, \Psi( X_s))
 -\<f, \Phi( X_s)\>\Big\}\d s ,\ \ f\in H^1.\end{split}\end{equation}

\beg{prp}\label{PP1}  $X_t$ is weakly c\'adl\'ag in $L^2(\mu)$,   $\eta_t$ is increasing and c\'adl\'ag in $H^{-1}$ and $(\ref{D1})$ holds.
\end{prp}
\beg{proof} (\ref{S'}) and  (\ref{DDE})   imply (\ref{D1}) and that $\eta_t$ is an increasing process in $H^{-1}$. In particular,
$$X_t= X_0 +\int_0^t \Big\{L\Psi(X_s)+\Phi(X_s)\Big\}\d s +\int_0^t\si_s\d W_s +\eta_t,\ \ t\ge 0$$ holds in $H^{-1}. $ Since the integral parts are continuous
in $H^{-1}$, it remains to show that $X_t$ is weakly c\'adl\'ag in $L^2(\mu)$ and hence c\'adl\'ag in $H^{-1}$ as $\sup_{t\in [0,T]}\|X_t\|^2<\infty$ for $T>0.$

Since $\bar\nu_t$ is increasing in $t$,
$$\bar\nu_{t+}:= \lim_{\vv\downarrow 0} \bar\nu_{t+\vv}\ge \bar\nu_t,\ \ t\ge 0.$$
Then, it is easy to see from (\ref {S'}) and {\bf (A1)}  that $X_t$ has weak left and right limits in $L^2(\mu)$ and its weak right limit $X_{t+}$ satisfies
\beq\label{CC}\<X_{t+}-  X_t, f\>= (\bar\nu_{t+}- \bar\nu_t)(f),\ \ \ f\in C_0(E).\end{equation}
Since $\bar\nu_{t+}\ge \bar\nu_t$, this in particular implies that $X_{t+}\ge X_t$.

On the other hand, by It\^o's formula and {\bf (A2)},
$$\|X_t^{(n)}\|^2-\|X_s^{(n)}\|^2\le 2\int_s^t \<\Phi(X_r^{(n)}), X_r^{(n)}\>\d r +\int_s^t \|\si_r\|_{HS}^2\d r + 2 \int_s^t \<X_r^{(n)}, \si_r\d W_r\>,\ \ 0\le s\le t.$$ Since $X^{(n)}\uparrow X$, by (\ref{3.2}), {\bf  (A2)}, {\bf (A3)},  and letting $n\uparrow\infty$, we obtain
$$\|X_t\|^2-\|X_s\|^2\le 2\int_s^t \<\Phi(X_r), X_r\>\d r +\int_s^t \|\si_r\|_{HS}^2\d r + 2 \int_s^t \<X_r, \si_r\d W_r\>,\ \ 0\le s\le t.$$ Therefore,
$$\|X_{t+}\|^2\le\liminf_{\vv\downarrow 0} \|X_{(t+\vv)\land T}\|^2 \le \|X_t\|^2.$$ Combining this with $X_{t+}\ge X_t$, we conclude that
$X_{t+}=X_t$, that is; $X_t$ is weakly right continuous in $L^2(\mu)$. \end{proof}

\subsection{Proof of Theorem \ref{T1.1}}

\paragraph{(a) Existence.} By Lemmas \ref{L3.1}, \ref{L3.2}, Proposition \ref{PP1} and (\ref{DDE}),  it remains to show that
$_{H^1}\< \Psi(X_t), \eta_t\>_{H^{-1}}=0, \d t$-a.e. Since $\Psi\in C^1$ with $\Psi(0)=0$ and $\Psi'\ge 0$, by $ X^{(n)} \uparrow  X\ge 0$ and (\ref{3.1}) we conclude  that (up to a subsequence)$$\Psi({X^{(n)}}^+)=\Psi(X^{(n)})^+\to \Psi(X)^+=\Psi(X)$$ weakly   $L^2([0,T]\to H^1;\d t).$ So,
\beg{equation*}\beg{split} &\int_0^T{_{H^1}\<\Psi(X_t),\eta_t\>_{H^{-1}}}\d t  =\lim_{n\to\infty}\int_0^T {_{H^1}\<\Psi((X_t^{(n)})^+), \eta_t\>_{H^{-1}}}\d t\\
 &=\lim_{n\to\infty} \lim_{m\to\infty} \int_0^T \Psi((X_t^{(n)})^+)(z) \bar\nu_t^{(m)}(\d z)\\
 & =\lim_{n\to\infty} \lim_{m\to\infty} m\int_0^T \Psi((X_t^{(n)})^+)
 (X_t^{(m)})^- \d t\d \mu =0,\ \ T>0.\end{split}\end{equation*} Since $_{H^1}\<\Psi(X),\eta\>_{H^{-1}}\ge 0$, we prove that $_{H^1}\< \Psi(X_t), \eta_t\>_{H^{-1}}=0, \d t$-a.e.

\paragraph{(b) The Markov property.}
For simplicity, we set $X_t= X_t(x)$. Let $0\le s_1<s_2<\cdots<s_m\le s<t$, and let $g\in C_b((H^{-1})^m)$. It remains to prove
\beq\label{MM1} \E\big\{f( X_t)g( X_{s_1},\cdots,  X_{s_m})\big\}= \E\big\{g( X_{s_1},\cdots,  X_{s_m})P_{t-s}f( X_s)\}\end{equation}
for any bounded and Lipschitz continuous $f$ in $L^1(\mu).$ By the Markov property of $X^{(n)}$ we have
$$ \E\big\{f(X_t^{(n)})g(X_{s_1}^{(n)},\cdots, X_{s_m}^{(n)})\big\}= \E\big\{g( X_{s_1}^{(n)},\cdots, X_{s_m}^{(n)})P_{t-s}^{(n)}f(X_{t-s}^{(n)})\},$$
where $P_{t-s}^{(n)}f(x):= \E f(X_t^{(n)}(x)).$
Since $X^{(n)}\uparrow  X$, and due to $\<\Phi(x)-\Phi(y)+nx^--ny^-, x-y\>\le l_0(x-y)^+$ and Theorem \ref{T2.3},
$P_{t-s}^{(n)}f$ is   continuous in $L^1(\mu)$ (hence also $L^2(\mu)$) uniformly w.r.t. $n\ge 1$, by letting $n\to \infty$ we obtain (\ref{MM1}).

\paragraph{(c) The $L^1$-Lipschitz continuity and consequences.}  Since (\ref{PH}) implies
$$\<\Phi(x)-\Phi(y)+nx^--ny^-, x-y\>\le K(x-y)^+,$$ by applying Theorem \ref{T2.3} to $X^{(n)}$ and letting $n\to\infty$, we prove (\ref{PH1}).
Then, for any Lipschitz continuous function $f$ on $L^1(\mu)$, and any $x\in L^1(\mu)$,
$$P_tf(x):=\lim_{n\to\infty} P_t f(x_n),\ \ \ x_n\to x\ \text{in}\ L^1(\mu),\ \text{and}\ \{x_n\}_{n\ge 1}\subset L^{1+r}(\mu)$$is well defined and provides
a Markov Lipschitz-Feller semigroup on $L^1(\mu).$  Moreover, by (\ref{LLL}) and It\^o's formula we have
$$\ff 1 T\int_0^T\|X_t(0)\|_{H^1}^2 \d t\le C,\ \ T>0,$$ for some constant $C>0$. Since $\|\cdot\|_{H^1}^2$ is a compact function in $L^1(\mu)$, this implies
that $P_t$ has an invariant probability measure $\pi$ with $\pi(\|\cdot\|_{H^1}^2)<\infty$. Finally, (\ref{PH2}) follows from (\ref{PH1}).

\paragraph{(d) Uniqueness.} Let $\Psi(s)= cs$ for some constant $c>0$, and let $(\tt X,\tt\eta)$ be another solution. We have
$$\d\|X_t-\tt X_t\|^2=2\big\{\<X_t-\tt X_t, \Phi(X_t)-\Phi(\tt X_t)\>-c \EE(X_t-\tt X_t, X_t-\tt X_t)\big\}\d t+ 2 _{H^1}\<X_t-\tt X_t, \d(\eta_t-\tt\eta_t)\>_{H^{-1}}.$$
Since $c>0$, $\Phi$ is Lipschitzian, $\d\eta_t,\d\tt\eta_t\ge 0$ and $$_{H^1}\<X_t,\d\eta_t\>_{H^{-1}}=_{H^1}\<\tt X_t,\d\tt\eta_t\>_{H^{-1}}=0,$$ this implies that
$$\d\|X_t-\tt X_t\|^2\le 2l_0 \|X_t-\tt X_t\|^2\d t.$$ Therefore, $X_t=\tt X_t$ holds for $t\in [0,T]$ provided $X_0=\tt X_0.$

\beg{thebibliography}{99}

\bibitem{0} V. Barbu, G. Da Prato, M. R\"ockner, \emph{Stochastic porous media equations and self-organized criticality,} Comm. Math. Phys. 285(2009), 901--923.

\bibitem{01} V. Barbu, G. Da Prato, M. R\"ockner, \emph{Finite time extinction for solutions to fast diffusion stochastic porous media equations,} C. R. Acad. Sci. Paris-Mathematics 347(2009), 81--84.

 \bibitem{02} V. Barbu, G. Da Prato, M. R\"ockner, \emph{Finite time extinction of solutions to fast diffusion equations driven by linear multiplicative noise,} J. Math. Anal. Appl. 389(2012), 147--164.

 \bibitem{03} A. Diaz-Guilera, \emph{Noise and dynamics of self-organized criticality phenomena,} Phys. Review A 45(1992), 8551--8558.

\bibitem{DP1} C. Donati-Martin, E. Pardoux, \emph{White noise driven SPDEs with reflection}. Probab. Theory Relat. Fields 95, 1-24(1993).

\bibitem{Gess} B. Gess, \emph{Strong solutions for stochastic partial
differential equations of gradient type,}  J. Funct. Anal. 263 (2012),   2355--2383.

\bibitem{GS} V. A. Gnatyuk, V. S. Shchirba, \emph{General properties of best approximation
with respect to a convex continuous function,} Ukrain. Mat. Zh. 34 (1982),
608--613.

\bibitem{MR} Z. Ma, M. R\"ockner, \emph{Introduction to the Theory of (Non-Symmetric) Dirichlet Forms,} Springer, Berlin, 1992.

\bibitem{NP}D. Nualart, E. Pardoux, \emph{White noise driven
quasilinear SPDEs with reflection}. Probab. Theory Relat. Fields
93,77-89(1992).

\bibitem{RRW07} J. Ren, M. R\"ockner, F.-Y. Wang, \emph{Stochastic generalized porous media and fast-diffusion equations, } J. Diff. Equations,  238(2007), 118--152.

\bibitem{RW08}	 M. R\"ockner, F.-Y. Wang, \emph{ Non-monotone stochastic generalized porous media equations,} J. Differential Equations 245(2008), 3898-3935.

\bibitem{7'}	 M. R\"ockner, F.-Y. Wang, \emph{General extinction results for stochastic partial differential equations and applications,} to appear in J. Lond. Math. Soc.

\bibitem{XZ} T. Xu and T.Zhang, \emph{White noise driven SPDEs with reflection:
existence, uniqueness and large deviation principles}. Stochastic Processes and Their Applications 119  (2009), 3453-3470.

\bibitem{ZA} L. Zambotti, \emph{A reflected stochastic heat equation as symmetric dynamics with respect to the 3-d Bessel bridge}. Journal of Functional Analysis 180, 195-209(2001).

\bibitem{Z} T.Zhang,  \emph{White noise driven SPDEs with reflection:
strong Feller properties and Harnack inequalities}. Potential Analysis 33 (2010), 137-151.

 \end{thebibliography}

\end{document}